# Sparse permutation invariant covariance estimation

**Adam J. Rothman**

*University of Michigan*
*Ann Arbor, MI 48109-1107*
*e-mail:* ajrothma@umich.edu

**Peter J. Bickel**

*University of California*
*Berkeley, CA 94720-3860*
*e-mail:* bickel@stat.berkeley.edu

**Elizaveta Levina**[*]

*University of Michigan*
*Ann Arbor, MI 48109-1107*
*e-mail:* elevina@umich.edu

**Ji Zhu**

*University of Michigan*
*Ann Arbor, MI 48109-1107*
*e-mail:* jizhu@umich.edu

**Abstract:** The paper proposes a method for constructing a sparse estimator for the inverse covariance (concentration) matrix in high-dimensional settings. The estimator uses a penalized normal likelihood approach and forces sparsity by using a lasso-type penalty. We establish a rate of convergence in the Frobenius norm as both data dimension $p$ and sample size $n$ are allowed to grow, and show that the rate depends explicitly on how sparse the true concentration matrix is. We also show that a correlation-based version of the method exhibits better rates in the operator norm. We also derive a fast iterative algorithm for computing the estimator, which relies on the popular Cholesky decomposition of the inverse but produces a permutation-invariant estimator. The method is compared to other estimators on simulated data and on a real data example of tumor tissue classification using gene expression data.

## 1. Introduction

Estimation of large covariance matrices, particularly in situations where the data dimension $p$ is comparable to or larger than the sample size $n$, has attracted a lot of attention recently. The abundance of high-dimensional data is one reason for the interest in the problem: gene arrays, fMRI, various kinds of spectroscopy, climate studies, and many other applications often generate very high dimensions and moderate sample sizes. Another reason is the ubiquity of the covariance matrix in data analysis tools. Principal component analysis (PCA), linear and quadratic discriminant analysis (LDA and QDA), inference about the means of the components, and analysis of independence and conditional independence in graphical models all require an estimate of the covariance matrix or its inverse, also known as the precision or concentration matrix. Finally, recent advances in random matrix theory – see Johnstone (2001) for a review, and also Paul (2007) – allowed in-depth theoretical studies of the traditional estimator, the sample (empirical) covariance matrix, and showed that without regularization the sample covariance performs poorly in high dimensions. These results helped stimulate research on alternative estimators in high dimensions.

Many alternatives to the sample covariance matrix have been proposed. A large class of methods covers the situation where variables have a natural ordering, e.g., longitudinal data, time series, spatial data, or spectroscopy. The implicit regularizing assumption underlying these methods is that variables far apart in the ordering have small correlations (or partial correlations, if the object of regularization is the concentration matrix). Methods for regularizing covariance by banding or tapering have been proposed by Bickel and Levina (2004) and Furrer and Bengtsson (2007). Bickel and Levina (2008) showed consistency of banded estimators in the operator norm under mild conditions as long as $(\log p)/n \to 0$, for both banding the covariance matrix and the Cholesky factor of the inverse discussed below.

When the inverse of the covariance matrix is the primary goal and the variables are ordered, regularization is usually introduced via the modified Cholesky decomposition,
$$\Sigma^{-1} = L^T D^{-1} L.$$
Here $L$ is a lower triangular matrix with $l_{jj} = 1$ and $l_{jj'} = -\phi_{jj'}$, where $\phi_{jj'}$, $j' < j$ is the coefficient of $X_{j'}$ in the population regression of $X_j$ on $X_1, \ldots, X_{j-1}$, and $D$ is a diagonal matrix with residual variances of these regressions on the diagonal. Several approaches to regularizing the Cholesky factor $L$ have been proposed, mostly based on its regression interpretation. A $k$-banded estimator of $L$ can be obtained by regressing each variable only on its closest $k$ predecessors; Wu and Pourahmadi (2003) proposed this estimator and chose $k$ via an AIC penalty. Bickel and Levina (2008) showed that banding the Cholesky factor produces a consistent estimator in the operator norm under weak conditions on the covariance matrix, and proposed a cross-validation scheme for picking $k$. Huang et al. (2006) proposed adding either an $l_2$ (ridge) or an $l_1$ (lasso) penalty on the elements of $L$ to the normal likelihood. The lasso penalty creates zeros in $L$ in arbitrary locations, which is more flexible than banding, but (unlike in the case of banding) the resulting estimate of the inverse may not have any zeros at all. Levina et al. (2008) proposed adaptive banding, which, by using a nested lasso penalty, allows a different $k$ for each regression, and hence is more flexible than banding while also retaining some sparsity in the inverse. Bayesian approaches to the problem introduce zeros via priors, either in the Cholesky factor (Smith and Kohn, 2002) or in the inverse itself (Wong et al., 2003).

There are, however, many applications where an ordering of the variables is not available: genetics, for example, or social and economic studies. Methods that are invariant to variable permutations (like the covariance matrix itself) are necessary in such applications. Regularizing large covariance matrices by Steinian shrinkage of eigenvalues has been proposed early on (Haff, 1980; Dey and Srinivasan, 1985). More recently, Ledoit and Wolf (2003) proposed a way to compute an optimal linear combination of the sample covariance with the identity matrix, which also results in shrinkage of eigenvalues. Shrinkage estimators are invariant to variable permutations but they do not affect the eigenvectors of the covariance, only the eigenvalues, and it has been shown that the sample eigenvectors are also not consistent when $p$ is large (Johnstone and Lu, 2004). Shrinking eigenvalues also does not create sparsity in any sense. Sometimes alternative estimators are available in the context of a specific application – e.g., for a factor analysis model with known factors Fan et al. (2008) develop regularized estimators for both the covariance and its inverse.

Our focus here will be on sparse estimators of the concentration matrix. Sparse concentration matrices are widely studied in the graphical models literature, since zero partial correlations imply a graph structure. The classical graphical models approach, however, is different from covariance estimation, since it normally focuses on just finding the zeros. For example, Drton and Perlman (2008) develop a multiple testing procedure for simultaneously testing hypotheses of zeros in the concentration matrix. There are also more algorithmic approaches to finding zeros in the concentration matrix, such as running a lasso regression of each variable on all the other variables (Meinshausen and Bühlmann, 2006), or the PC-algorithm (Kalisch and Bühlmann, 2007). Both have been shown to be consistent in high-dimensional settings, but none of these methods supply an estimator of the covariance matrix. In principle, once the zeros are found, a constrained maximum likelihood estimator of the covariance can be computed (Chaudhuri et al., 2007), but it is not clear what the properties of such a two-step procedure would be.

Two recent papers, d'Aspremont et al. (2008) and Yuan and Lin (2007), take a penalized likelihood approach by applying an $l_1$ penalty to the entries of the concentration matrix. This results in a permutation-invariant loss function that tends to produce a sparse estimate of the inverse. Yuan and Lin (2007) used the max-det algorithm to compute the estimator, which limited their numerical results to values of $p \leq 10$, and derived a fixed $p$, large $n$ convergence result. d'Aspremont et al. (2008) proposed a much faster semi-definite programming algorithm based on Nesterov's method for interior point optimization. While this paper was in review, a new very fast algorithm for the same problem was proposed by Friedman et al. (2008), which is based on the coordinate descent algorithm for the lasso (Friedman et al., 2007).

In this paper, we analyze the estimator resulting from penalizing the normal likelihood with the $l_1$ penalty on the entries of the concentration matrix (we will refer to this estimator as SPICE – Sparse Permutation Invariant Covariance Estimator) in the high-dimensional setting, allowing both the dimension $p$ and the sample size $n$ to grow. We give an explicit convergence rate in the Frobenius norm and show that the rate depends on how sparse the true concentration matrix is. For a slight modification of the method based on using the sample correlation matrix, we obtain the rate of convergence in operator norm and show that it is essentially equivalent to the rate of thresholding the covariance matrix itself obtained in Bickel and Levina (2007). We also derive our own optimization algorithm for computing the estimator, based on Cholesky decomposition and the local quadratic approximation. Unlike other estimation methods that rely on the Cholesky decomposition, our algorithm is invariant under variable permutations. Because we use the local quadratic approximation, the algorithm is equally applicable to general $l_q$ penalties on the entries of the inverse, not just $l_1$.

The rest of the paper is organized as follows: Section 2 summarizes the SPICE approach in general, and presents consistency results. The Cholesky-based computational algorithm, along with a discussion of optimization issues, is presented in Section 3. Section 4 presents numerical results for SPICE and a number of other methods, for simulated data and a real example on classification of colon tumors using gene expression data. Section 5 concludes with discussion.

## 2. Analysis of the SPICE method

We assume throughout that we observe $\boldsymbol{X}_1, \ldots, \boldsymbol{X}_n$, i.i.d. $p$-variate normal random variables with mean $\boldsymbol{0}$ and covariance matrix $\Sigma_0$, and write $\boldsymbol{X}_i = (X_{i1}, \ldots, X_{ip})^T$. Let $\Sigma_0 = [\sigma_{0ij}]$, and $\Omega_0 = \Sigma_0^{-1}$ be the inverse of the true covariance matrix. For any matrix $M = [m_{ij}]$, we write $|M|$ for the determinant of $M$, $\text{tr}(M)$ for the trace of $M$, and $\varphi_{\max}(M)$ and $\varphi_{\min}(M)$ for the largest and smallest eigenvalues, respectively. We write $M^+ = \text{diag}(M)$ for a diagonal matrix with the same diagonal as M, and $M^- = M - M^+$. In the asymptotic analysis, we will use the Frobenius matrix norm $\|M\|_F^2 = \sum_{i,j} m_{ij}^2$, and the operator norm (also known as matrix 2-norm), $\|M\|^2 = \varphi_{\max}(MM^T)$. We will also write $|\cdot|_1$ for the $l_1$ norm of a vector or matrix vectorized, i.e., for a matrix $|M|_1 = \sum_{i,j} |m_{ij}|$.

It is easy to see that under the normal assumption the negative log-likelihood, up to a constant, can be written in terms of the concentration matrix as

$$\ell(\boldsymbol{X}_1, \ldots, \boldsymbol{X}_n; \Omega) = \text{tr}(\Omega \hat{\Sigma}) - \log |\Omega|,$$

where
$$\hat{\Sigma} = \frac{1}{n} \sum_{i=1}^{n} (\boldsymbol{X}_i - \bar{\boldsymbol{X}})(\boldsymbol{X}_i - \bar{\boldsymbol{X}})^T$$

is the sample covariance matrix.

We define the SPICE estimator $\hat{\Omega}_\lambda$ of the inverse covariance matrix as the minimizer of the penalized negative log-likelihood,

$$\hat{\Omega}_\lambda = \arg\min_{\Omega \succ 0} \left\{ \text{tr}(\Omega \hat{\Sigma}) - \log|\Omega| + \lambda |\Omega^-|_1 \right\} \tag{1}$$

where $\lambda$ is a non-negative tuning parameter, and the minimization is taken over symmetric positive definite matrices.

SPICE is identical to the lasso-type estimator proposed by Yuan and Lin (2007), and very similar to the estimator of d'Aspremont et al. (2008) (they used $|\Omega|_1$ rather than $|\Omega^-|_1$ in the penalty). The loss function is invariant to permutations of variables and should encourage sparsity in $\hat{\Omega}$ due to the $l_1$ penalty applied to its off-diagonal elements.

We make the following assumptions about the true model:

A1: Let the set $S = \{(i, j) : \Omega_{0ij} \neq 0, i \neq j\}$. Then $\text{card}(S) \leq s$.
A2: $\varphi_{\min}(\Sigma_0) \geq \underline{k} > 0$, or equivalently $\varphi_{\max}(\Omega_0) \leq 1/\underline{k}$.
A3: $\varphi_{\max}(\Sigma_0) \leq \overline{k}$.

Note that assumption A2 guarantees that $\Omega_0$ exists. Assumption A1 is more of a definition, since it does not stipulate anything about $s$ ($s = p(p-1)/2$ would give a full matrix).

**Theorem 1.** *Let $\hat{\Omega}_\lambda$ be the minimizer defined by (1). Under A1, A2, A3, if $\lambda \asymp \sqrt{\frac{\log p}{n}}$,*

$$\|\hat{\Omega}_\lambda - \Omega_0\|_F = O_P\left(\sqrt{\frac{(p+s)\log p}{n}}\right). \tag{2}$$

The theorem can be restated, more suggestively, as

$$\frac{\|\hat{\Omega}_\lambda - \Omega_0\|_F^2}{p} = O_P\left(\left(1 + \frac{s}{p}\right)\frac{\log p}{n}\right). \tag{3}$$

The reason for the second formulation (3) is the relation of the Frobenius norm to the operator norm, $\|M\|_F^2/p \leq \|M\|^2 \leq \|M\|_F^2$.

Before proceeding with the proof of Theorem 1, we discuss a modification to SPICE based on using the correlation matrix. An inspection of the proof reveals that the worst part of the rate, $\sqrt{p \log p/n}$, comes from estimating the diagonal. This suggests that if we were to use the correlation matrix rather than the covariance matrix, we should be able to get the rate of $\sqrt{s \log p/n}$. Indeed, let $\Sigma_0 = W\Gamma W$, where $\Gamma$ is the true correlation matrix, and $W$ is the diagonal matrix of true standard deviations. Let $\hat{W}$ and $\hat{\Gamma}$ be the sample estimates of $W$ and $\Gamma$, i.e., $\hat{W}^2 = \hat{\Sigma}^+$, $\hat{\Gamma} = \hat{W}^{-1}\hat{\Sigma}\hat{W}^{-1}$. Let $K = \Gamma^{-1}$. Define a SPICE estimate of $K$ by

$$\hat{K}_\lambda = \arg\min_{\Omega \succ 0} \left\{ \mathrm{tr}(\Omega\hat{\Gamma}) - \log|\Omega| + \lambda|\Omega^-|_1 \right\} \quad (4)$$

Then we can define a modified correlation-based estimator of the concentration matrix by

$$\tilde{\Omega}_\lambda = \hat{W}^{-1}\hat{K}_\lambda\hat{W}^{-1}. \quad (5)$$

It turns out that in the Frobenius norm $\tilde{\Omega}$ has the same rate as $\hat{\Omega}$, but for $\tilde{\Omega}$ we can get a convergence rate in the operator norm (matrix 2-norm). As discussed previously by Bickel and Levina (2008), El Karoui (2007) and others, the operator norm is more appropriate than the Frobenius norm for spectral analysis, e.g., PCA. It also allows for a direct comparison with banding rates obtained in Bickel and Levina (2008) and thresholding rates in Bickel and Levina (2007).

**Theorem 2.** *Under assumptions of Theorem 1,*

$$\|\tilde{\Omega}_\lambda - \Omega_0\| = O_P\left(\sqrt{\frac{(s+1)\log p}{n}}\right).$$

*Note.* This rate is very similar to the rate for thresholding the covariance matrix obtained by Bickel and Levina (2007). They showed that under the assumption $\max_i \sum_j |\sigma_{ij}|^q \leq c_0(p)$ for $0 \leq q < 1$, if the sample covariance entries are set to 0 when their absolute values fall below the threshold $\lambda = M\sqrt{\frac{\log p}{n}}$, then the resulting estimator converges to the truth in operator norm at the rate no worse than $O_P\left(c_0(p)(\frac{\log p}{n})^{(1-q)/2}\right)$. Since the truly sparse case corresponds to $q = 0$, and $c_0(p)$ is a bound on the number of non-zero elements in each row, and thus $\sqrt{s} \asymp c_0(p)$, this rate coincides with ours, even though the estimator and the method of proof are very different. However, Lemma 1 below is the basis of the proof in both cases, and ultimately it is the bound (6) that gives rise to the same rate. A similar rate has been obtained for banding the covariance matrix in Bickel and Levina (2008), under an additional assumption that depends on the ordering of the variables and is not applicable here (see Bickel and Levina (2007) for a comparison between banding and thresholding rates).

In the proof, we will need a lemma of Bickel and Levina (2008) (Lemma 3) which is based on a large deviation result of Saulis and Statulevičius (1991). We state the result here for completeness.

**Lemma 1.** *Let $Z_i$ be i.i.d. $\mathcal{N}(\mathbf{0}, \Sigma_p)$ and $\varphi_{\max}(\Sigma_p) \leq \overline{k} < \infty$. Then, if $\Sigma_p = [\sigma_{ab}]$,*

$$P\left[\left|\sum_{i=1}^n (Z_{ij}Z_{ik} - \sigma_{jk})\right| \geq n\nu\right] \leq c_1 \exp(-c_2 n\nu^2) \quad \textit{for } |\nu| \leq \delta \quad (6)$$

*where $c_1$, $c_2$ and $\delta$ depend on $\overline{k}$ only.*

*Proof of Theorem 1.* Let

$$Q(\Omega) = \text{tr}(\Omega\hat{\Sigma}) - \log|\Omega| + \lambda|\Omega^-|_1 - \text{tr}(\Omega_0\hat{\Sigma}) + \log|\Omega_0| - \lambda|\Omega_0^-|_1$$
$$= \text{tr}\big[(\Omega - \Omega_0)(\hat{\Sigma} - \Sigma_0)\big] - (\log|\Omega| - \log|\Omega_0|)$$
$$+ \text{tr}\big[(\Omega - \Omega_0)\Sigma_0\big] + \lambda(|\Omega^-|_1 - |\Omega_0^-|_1) \quad (7)$$

Our estimate $\hat{\Omega}$ minimizes $Q(\Omega)$, or equivalently $\hat{\Delta} = \hat{\Omega} - \Omega_0$ minimizes $G(\Delta) \equiv Q(\Omega_0 + \Delta)$. Note that we suppress the dependence on $\lambda$ in $\hat{\Omega}$ and $\hat{\Delta}$.

The main idea of the proof is as follows. Consider the set

$$\Theta_n(M) = \{\Delta : \Delta = \Delta^T, \|\Delta\|_F = Mr_n\},$$

where

$$r_n = \sqrt{\frac{(p+s)\log p}{n}} \to 0.$$

Note that $G(\Delta) = Q(\Omega_0 + \Delta)$ is a convex function, and

$$G(\hat{\Delta}) \leq G(0) = 0.$$

Then, if we can show that

$$\inf\{G(\Delta) : \Delta \in \Theta_n(M)\} > 0,$$

the minimizer $\hat{\Delta}$ must be inside the sphere defined by $\Theta_n(M)$, and hence

$$\|\hat{\Delta}\|_F \leq Mr_n. \quad (8)$$

For the logarithm term in (7), doing the Taylor expansion of $f(t) = \log|\Omega + t\Delta|$ and using the integral form of the remainder and the symmetry of $\Delta$, $\Sigma_0$, and $\Omega_0$ gives

$$\log|\Omega_0 + \Delta| - \log|\Omega_0| = \text{tr}(\Sigma_0\Delta) - \tilde{\Delta}^T\left[\int_0^1 (1-v)(\Omega_0 + v\Delta)^{-1} \otimes (\Omega_0 + v\Delta)^{-1} dv\right]\tilde{\Delta} \quad (9)$$

where $\otimes$ is the Kronecker product (if $A = [a_{ij}]_{p_1 \times q_1}$, $B = [b_{kl}]_{p_2 \times q_2}$, then $A \otimes B = [a_{ij}b_{kl}]_{p_1p_2 \times q_1q_2}$), and $\tilde{\Delta}$ is $\Delta$ vectorized to match the dimensions of the Kronecker product.

Therefore, we may write (7) as,

$$G(\Delta) = \text{tr}\big(\Delta(\hat{\Sigma} - \Sigma_0)\big) + \tilde{\Delta}^T\left[\int_0^1 (1-v)(\Omega_0 + v\Delta)^{-1} \otimes (\Omega_0 + v\Delta)^{-1} dv\right]\tilde{\Delta}$$
$$+ \lambda(|\Omega_0^- + \Delta^-|_1 - |\Omega_0^-|_1) \quad (10)$$

For an index set $A$ and a matrix $M = [m_{ij}]$, write $M_A \equiv [m_{ij}I((i,j) \in A)]$, where $I(\cdot)$ is an indicator function. Recall $S = \{(i,j) : \Omega_{0ij} \neq 0, i \neq j\}$ and let $\overline{S}$ be its complement. Note that $|\Omega_0^- + \Delta^-|_1 = |\Omega_{0S}^- + \Delta_S^-|_1 + |\Delta_{\overline{S}}^-|_1$, and $|\Omega_0^-|_1 = |\Omega_{0S}^-|_1$. Then the triangular inequality implies

$$\lambda(|\Omega_0^- + \Delta^-|_1 - |\Omega_0^-|_1) \geq \lambda(|\Delta_{\overline{S}}^-|_1 - |\Delta_S^-|_1). \quad (11)$$

Now, using symmetry again, we write

$$|\text{tr}(\Delta(\hat{\Sigma} - \Sigma_0))| \leq \left|\sum_{i \neq j}(\hat{\sigma}_{ij} - \sigma_{0ij})\Delta_{ij}\right| + \left|\sum_{i}(\hat{\sigma}_{ii} - \sigma_{0ii})\Delta_{ii}\right| = \text{I} + \text{II}. \quad (12)$$

To bound term I, note that the union sum inequality and Lemma 1 imply that, with probability tending to 1,

$$\max_{i \neq j}|\hat{\sigma}_{ij} - \sigma_{0ij}| \leq C_1\sqrt{\frac{\log p}{n}}$$

and hence term I is bounded by

$$\text{I} \leq C_1\sqrt{\frac{\log p}{n}}|\Delta^-|_1. \quad (13)$$

The second bound comes from the Cauchy-Schwartz inequality and Lemma 1:

$$\text{II} \leq \left[\sum_{i=1}^{p}(\hat{\sigma}_{ii} - \sigma_{0ii})^2\right]^{1/2}\|\Delta^+\|_F \leq \sqrt{p}\max_{1\leq i \leq p}|\hat{\sigma}_{ii} - \sigma_{0ii}|\|\Delta^+\|_F$$

$$\leq C_2\sqrt{\frac{p\log p}{n}}\|\Delta^+\|_F \leq C_2\sqrt{\frac{(p+s)\log p}{n}}\|\Delta^+\|_F, \quad (14)$$

also with probability tending to 1.

Now, take

$$\lambda = \frac{C_1}{\varepsilon}\sqrt{\frac{\log p}{n}}. \quad (15)$$

By (10),

$$G(\Delta) \geq \frac{1}{4}\underline{k}^2\|\Delta\|_F^2 - C_1\sqrt{\frac{\log p}{n}}|\Delta^-|_1 - C_2\sqrt{\frac{(p+s)\log p}{n}}\|\Delta^+\|_F$$
$$+ \lambda(|\Delta_{\overline{S}}^-|_1 - |\Delta_S^-|_1)$$
$$= \frac{1}{4}\underline{k}^2\|\Delta\|_F^2 - C_1\sqrt{\frac{\log p}{n}}\left(1 - \frac{1}{\varepsilon}\right)|\Delta_{\overline{S}}^-|_1 - C_1\sqrt{\frac{\log p}{n}}\left(1 + \frac{1}{\varepsilon}\right)|\Delta_S^-|_1$$
$$- C_2\sqrt{\frac{(p+s)\log p}{n}}\|\Delta^+\|_F \quad (16)$$

The first term comes from a bound on the integral which we will argue separately below. The second term is always positive, and hence we may omit it for the lower bound. Now, note that

$$|\Delta_S^-|_1 \leq \sqrt{s}\|\Delta_S^-\|_F \leq \sqrt{s}\|\Delta^-\|_F \leq \sqrt{p+s}\|\Delta^-\|_F.$$

Thus we have

$$\begin{aligned} G(\Delta) &\geq \|\Delta^-\|_F^2 \left[\frac{1}{4}\underline{k}^2 - C_1\sqrt{\frac{(p+s)\log p}{n}}\left(1 + \frac{1}{\varepsilon}\right)\|\Delta^-\|_F^{-1}\right] \\ &\quad + \|\Delta^+\|_F^2 \left[\frac{1}{4}\underline{k}^2 - C_2\sqrt{\frac{(p+s)\log p}{n}}\|\Delta^+\|_F^{-1}\right] \\ &= \|\Delta^-\|_F^2 \left[\frac{1}{4}\underline{k}^2 - \frac{C_1(1+\varepsilon)}{\varepsilon M}\right] + \|\Delta^+\|_F^2 \left[\frac{1}{4}\underline{k}^2 - \frac{C_2}{M}\right] > 0 \quad (17) \end{aligned}$$

for $M$ sufficiently large.

It only remains to check the bound on the integral term in (10). Recall that $\varphi_{\min}(M) = \min_{\|x\|=1} x^T M x$. After factoring out the norm of $\tilde{\Delta}$, we have, for $\Delta \in \Theta_n(M)$,

$$\begin{aligned} \varphi_{\min}&\left(\int_0^1 (1-v)(\Omega_0 + v\Delta)^{-1} \otimes (\Omega_0 + v\Delta)^{-1} dv\right) \\ &\geq \int_0^1 (1-v)\varphi_{\min}^2(\Omega_0 + v\Delta)^{-1} dv \geq \frac{1}{2} \min_{0 \leq v \leq 1} \varphi_{\min}^2(\Omega_0 + v\Delta)^{-1} \\ &\geq \frac{1}{2} \min\left\{\varphi_{\min}^2(\Omega_0 + \Delta)^{-1} : \|\Delta\|_F \leq Mr_n\right\}. \end{aligned}$$

The first inequality uses the fact that the eigenvalues of Kronecker products of symmetric matrices are the products of the eigenvalues of their factors. Now

$$\varphi_{\min}^2(\Omega_0 + \Delta)^{-1} = \varphi_{\max}^{-2}(\Omega_0 + \Delta) \geq (\|\Omega_0\| + \|\Delta\|)^{-2} \geq \frac{1}{2}\underline{k}^2 \quad (18)$$

with probability tending to 1, since $\|\Delta\| \leq \|\Delta\|_F = o(1)$. This establishes the theorem. □

As noted above, an inspection of the proof shows that $\sqrt{p \log p/n}$ in the rate comes from estimating the diagonal. If we focus on the correlation matrix estimate $\hat{K}_\lambda$ in (4) instead, we can immediately obtain

**Corollary 1.** *Under assumptions of Theorem 1,*

$$\|\hat{K}_\lambda - K\|_F = O_P\left(\sqrt{\frac{s \log p}{n}}\right).$$

Now we can use Corollary 1 to prove Theorem 2, the operator norm bound.

*Proof of Theorem 2.* Write

$$\begin{aligned} \|\tilde{\Omega}_\lambda - \Omega_0\| &= \|\hat{W}^{-1}\hat{K}_\lambda \hat{W}^{-1} - W^{-1}KW^{-1}\| \\ &\leq \|\hat{W}^{-1} - W^{-1}\|\,\|\hat{K}_\lambda - K\|\,\|\hat{W}^{-1} - W^{-1}\| \\ &\quad + \|\hat{W}^{-1} - W^{-1}\|(\|\hat{K}_\lambda\|\,\|W^{-1}\| + \|\hat{W}^{-1}\|\,\|K\|) \\ &\quad + \|\hat{K}_\lambda - K\|\,\|\hat{W}^{-1}\|\,\|W^{-1}\| \end{aligned}$$

where we are using the sub-multiplicative norm property $\|AB\| \leq \|A\|\|B\|$ (see, e.g., Golub and Van Loan (1989)). Now, $\|W^{-1}\|$ and $\|K\|$ are $O(1)$ by assumptions A2 and A3. Lemma 1 implies that

$$\|\hat{W}^2 - W^2\| = O_P\left(\sqrt{\frac{\log p}{n}}\right), \qquad (19)$$

and since $\|\hat{W}^{-1} - W^{-1}\| \stackrel{P}{\asymp} \|\hat{W}^2 - W^2\|$ (where by $A \stackrel{P}{\asymp} B$ we mean $A = O_P(B)$ and $B = O_P(A)$), we have the rate of $\sqrt{\log p/n}$ for $\|\hat{W}^{-1} - W^{-1}\|$. This together with Corollary 1 in turn implies that $\|\hat{W}^{-1}\|$ and $\|\hat{K}_\lambda\|$ are $O_P(1)$, and the theorem follows. □

Note that in the Frobenius norm, we only have $\|\hat{W}^2 - W^2\| = O_P(\sqrt{p \log p/n})$, and thus the Frobenius rate of $\tilde{\Omega}_\lambda$ is the same as that of $\hat{\Omega}_\lambda$.

## 3. The Cholesky-based SPICE algorithm

In this section, we develop an iterative algorithm for computing the SPICE estimator using the Cholesky decomposition; however, unlike other estimators that depend on the Cholesky decomposition, we minimize a permutation invariant objective function, and thus the estimator remains permutation invariant. We use the quadratic approximation to the absolute value, a standard tool in optimization which has been previously used in the statistics literature to handle lasso-type penalties, for example, by Fan and Li (2001) and Huang et al. (2006). In this our algorithm differs from the glasso algorithm of Friedman et al. (2008), which is based on a lasso algorithm and works directly on the absolute values. Both algorithms have computation complexity of $O(p^3)$, but we acquire another small constant factor (on the order of 10) due to the additional iterations required for the quadratic approximation to converge (see more on this in Section 4). However, using the quadratic approximation allows us to write down the algorithm explicitly in general terms for an $l_q$ penalty $|w_{ij}|^q$ with $q \geq 1$, rather than only for $q = 1$. In particular, our algorithm is equally applicable for use with a ridge penalty ($q = 2$), although in that special case it simplifies even further, or with a bridge penalty ($1 < q < 2$) proposed by Fu (1998), which may work better for certain classes of covariances. It can also be used with SCAD (Fan and Li, 2001) or other more complicated non-convex penalties that are typically approximated by the local quadratic approximation. Even though we derive the algorithm with a general $q$, in this paper we only present results for $q = 1$.

Our goal is to minimize the objective function,

$$f(\Omega) = \text{tr}(\Omega \hat{\Sigma}) - \log |\Omega| + \lambda \sum_{j' \neq j} |\omega_{j'j}|^q, \qquad (20)$$

where $q = 1$ corresponds to the computation of $\hat{\Omega}_\lambda$ in (1). For $q \geq 1$, the objective function is convex in the elements of $\Omega$ and has a global minimum $\hat{\Omega}$.

*A.J. Rothman et al./Sparse covariance estimation*      504
Our strategy is to re-parametrize the objective (20) using the Cholesky decomposition of $\Omega$ to enforce automatic positive definiteness. Rather than using the modified Cholesky decomposition with its regression interpretation, as has been standard in the literature, we simply write

$$\Omega = T^T T,$$

where $T = [t_{ij}]$ is a lower triangular matrix. We can still use the regression interpretation if needed, by writing

$$\begin{aligned} t_{jj'} &= -\frac{\phi_{jj'}}{\sqrt{d_{jj}}}, \quad j' < j \\ t_{jj} &= \frac{1}{\sqrt{d_{jj}}}, \end{aligned} \quad (21)$$

where $\phi_{jj'}$ is the coefficient of $X_{j'}$ in the regression of $X_j$ on $X_1, \ldots, X_{j-1}$, and $d_{jj}$ is the corresponding residual variance.

To minimize $f$ in terms of $T$, we apply a cyclical coordinate descent approach and minimize $f$ with respect to one element of $T$ at a time. Further, we use a quadratic approximation to $f$, which allows us to find the minimum of the univariate functions of each parameter in closed form. The algorithm is iterated until convergence. Here we outline the main steps of the algorithm, and leave the full derivation for the Appendix.

In a slight abuse of notation, we write $X$ for the $n \times p$ data matrix where each column has already been centered by its sample mean. The three terms in (20) can be expressed as a function of $T$ as follows:

$$\mathrm{tr}(\Omega\hat{\Sigma}) = \frac{1}{n} \sum_{i=1}^{n} \sum_{j=1}^{p} \left( \sum_{k=1}^{j} t_{jk} X_{ik} \right)^2 \quad (22)$$

$$\log|\Omega| = 2 \sum_{j=1}^{p} \log t_{jj} \quad (23)$$

$$\sum_{j' \neq j} |\omega_{j'j}|^q = 2 \sum_{j' > j} \left| \sum_{k=j'}^{p} t_{kj'} t_{kj} \right|^q \quad (24)$$

The quadratic approximation for $|u|^q$ is shown in (25). Since the algorithm is iterative, $u^{(k)}$ denotes the value of $u$ from the previous iteration, and $u^{(k+1)}$ is the value at current iteration.

$$|u^{(k+1)}|^q \approx \frac{q}{2} \frac{(u^{(k+1)})^2}{|u^{(k)}|^{2-q}} + \left(1 - \frac{q}{2}\right) |u^{(k)}|^q \quad (25)$$

Hunter and Li (2005) suggest replacing $|u^{(k)}|$ in the denominator with $|u^{(k)}| + \epsilon$ to avoid division by zero, and refer to this as the $\epsilon$-perturbed quadratic approximation. This quadratic approximation to $f$, which we denote $\tilde{f}_{\epsilon,k}$ at iteration $k$, allows us to easily take partial derivatives with respect to each parameter in $T$, and provides a closed form solution for the univariate minimizer for each coordinate.

The algorithm requires an initial value $\hat{T}^{(0)}$, which corresponds to $\hat{\Omega}^{(0)}$. If the sample covariance $\hat{\Sigma}$ is non-degenerate, which is generally the case for $p < n$, one could simply set $\hat{\Omega}^{(0)} = \hat{\Sigma}^{-1}$. More generally, we found the following simple strategy to work well: approximate $\phi_{jj'}$ in (21) by regressing $X_j$ on $X_{j'}$ *alone*, for $j' = 1, \ldots, j-1$, and then compute $\hat{T}^{(0)}$ using (21). Yet another alternative is to start from the diagonal estimator.

**The Algorithm:**

Step 0. Initialize $\hat{T} = \hat{T}^{(0)}$ and $\hat{\Omega}^{(0)} = (\hat{T}^{(0)})^T \hat{T}^{(0)}$.

Step 1. For each parameter $t_{lc}$, $c = 1, \ldots, p, l = c, \ldots, p$, solve $\nabla_{t_{lc}} \tilde{f}_{\epsilon,k}(T) = 0$ to find new $\hat{t}_{lc}$.

Step 2. Repeat Step 1 until convergence of $\hat{T}$ and set $T^{(k+1)} = \hat{T}$.

Step 3. Set $\hat{\Omega}^{(k+1)} = (T^{(k+1)})^T T^{(k+1)}$ and repeat Steps 1–3 until convergence of $\hat{\Omega}$.

Steps 2 and 3 may seem redundant, but they are needed for two different reasons. Step 2 is needed because we only minimize with respect to one parameter at a time, holding all other parameters fixed; and Step 3 is needed because of the quadratic approximation for $|u|^q$. After convergence, we replace entries in $\hat{\Omega}$ with smaller magnitude than $\epsilon$ with zero, using a fixed value of $\epsilon = 10^{-8}$. Another approach with virtually the same performance is to replace entries of $\hat{\Omega}^{(k)}$ with $\epsilon$ if their magnitude falls below $\epsilon$ in Step 3, and use (25) directly in the objective function in Step 1 instead of using $\tilde{f}_{\epsilon,k}$.

In practice, we found that working with the correlation matrix as described in Theorem 2 is slightly better than working with the covariance matrix, although the differences are fairly small. Still, in all the numerical results we standardize the variables first and then rescale our estimate by the sample standard deviations of the variables.

### *3.1. Algorithm convergence*

The convergence of the algorithm essentially follows from two standard results. For the inner loop cycling through individual parameters, the value of the objective function decreases at each iteration, and the objective function is differentiable everywhere. Thus the inner loop of the algorithm converges by a standard theorem on cyclical coordinate descent for smooth functions (see, e.g., Bazaraa et al. (2006), p. 367), to a stationary point $\nabla g(T) = 0$, where $g(T) = \tilde{f}_{\epsilon,k}(T^T T)$. The function $\tilde{f}_{\epsilon,k}$ is convex in the original parameters $\omega_{ij}$, but since we reparametrized it in terms of $T$, the function $g$ is not necessarily convex in $T$. In the next proposition we verify that this stationary point of $g$ corresponds to the global minimum of the convex function $\tilde{f}_{\epsilon,k}$.

**Proposition 1.** *Let $\tilde{f} \equiv \tilde{f}_{\epsilon,k}$ be the original convex function $f$ approximated by the $\epsilon$-perturbed local quadratic approximation at iteration $k$, let $T$ be a $p \times p$* *lower triangular matrix, and let* $g(T) = \tilde{f}(T^T T)$. *Let* $S_0$ *be the unique solution to* $\nabla \tilde{f}(S) = 0$, *and let* $T_0$ *be a solution to* $\nabla g(T) = 0$. *Then* $S_0 = T_0^T T_0$.

*Proof of Proposition 1.* Let $h : T \to T^T T$. Note that $h$ maps all of $\mathbb{R}^{p(p+1)/2}$ (all lower triangular matrices) into a convex subset of $\mathbb{R}^{p(p+1)/2}$ (non-negative definite symmetric matrices). Denote the differential of $h$ in the direction $d \in \mathbb{R}^{p(p+1)/2}$ evaluated at $t_0 \in \mathbb{R}^{p(p+1)/2}$ by $\nabla h(t_0)[d]$. Then,

$$\nabla h(t_0)[d] = T_0^T D + D^T T_0, \tag{26}$$

where $T_0$ and $D$ are, respectively, $t_0$ and $d$ written as $p \times p$ matrices. Now, using the chain rule and (26), we have

$$\nabla g(t_0)[d] = \nabla \tilde{f}(vec(T_0^T T_0)) \left( T_0^T D + D^T T_0 \right). \tag{27}$$

where we now think of $\tilde{f}$ as a function from $\mathbb{R}^{p(p+1)/2}$ to $\mathbb{R}$. Since $\tilde{f}$ is convex and has a unique minimizer $s_0 = vec(S_0)$, $\nabla \tilde{f}(s)[d]$ vanishes iff $s = s_0$ or $d = 0$. Thus $\nabla g(t_0)[d] = 0$ vanishes iff $T_0^T T_0 = S_0$ or $T_0^T D + D^T T_0 = 0$, or $T_0^T D = -(T_0^T D)^T$. If any diagonal elements of $T_0$ are 0, then $T_0$ is singular, and so is $T_0^T T_0$, and thus $g(T_0) = \infty$, so a singular $T_0$ cannot be a stationary point of $g$. Since $T_0$ is lower triangular and all its diagonal elements must be non-zero, one can show by induction that $T_0^T D = -(T_0^T D)^T$ implies $D = 0$. □

For the outer loop iterating through the quadratic approximation, we can apply the argument of Hunter and Li (2005) for $\epsilon$-perturbed local quadratic approximation obtained from general results for minorize-maximize algorithms, and conclude that as $k \to \infty$ and $\epsilon \to 0$ the algorithm converges to the global minimum of the original convex function $f$ in (20). In practice, we have also observed that our algorithm and glasso converge to the same solution.

### *3.2. Computational complexity*

The computational complexity of the algorithm in terms of $p$ is $O(p^3)$, since each parameter update is at most $O(p)$ (see (32) in the Appendix), and there are $O(p^2)$ parameters. The only other algorithm for computing this estimator at the cost of $O(p^3)$ is glasso of Friedman et al. (2008); the algorithms of Yuan and Lin (2007) and d'Aspremont et al. (2008) have higher computational cost. For extensive timing comparisons of glasso and the algorithm of d'Aspremont et al. (2008), which showed convincingly that glasso is much faster, see Friedman et al. (2008). The exact timing also depends on the implementation, platform, etc (our algorithm is implemented in C and glasso in Fortran). Actual computing times we obtained for glasso and the SPICE algorithm are shown below in Figure 1, for model $\Omega_2$ described in Section 4.1, with values of tuning parameters chosen as described in Section 3.3.

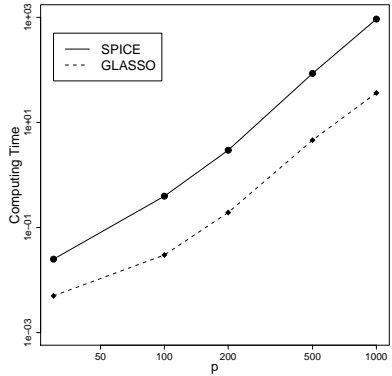

Fig 1. *Computing time in seconds vs p (log-log scale) for SPICE and glasso.*

### *3.3. Choice of tuning parameter*

Like any other penalty-based approach, SPICE requires selecting the tuning parameter $\lambda$. In simulations, we generate a separate validation dataset, and select $\lambda$ by maximizing the normal likelihood on the validation data with $\hat{\Omega}_\lambda$ estimated from the training data. Alternatively, one can use 5-fold cross-validation, which we do for the real data analysis. There is some theoretical basis for selecting the tuning parameter in this way – see Bickel and Levina (2007).

## 4. Numerical Results

In this section, we compare the performance of SPICE to the shrinkage estimator of Ledoit and Wolf (2003) and to the sample covariance matrix when applicable ($p < n$), using simulated and real data. We do not include any estimators that depend on variable ordering (such as banding of Bickel and Levina (2008) or the Lasso penalty on the Cholesky factor of Huang et al. (2006)), nor estimators that focus on introducing sparsity in the covariance matrix itself rather than in its inverse (such as thresholding), as they would automatically be at a disadvantage on sparse concentration matrices. The Ledoit-Wolf estimator does not introduce sparsity in the inverse either, but we use it as a benchmark for cases when $p > n$, since the sample covariance is not invertible.

### *4.1. Simulations*

In simulations, we focus on comparing performance on sparse concentration matrices, with varying levels of sparsity. We consider the following four covariance models.

1. $\Omega_1$: AR(1), $\sigma_{j'j} = 0.7^{|j'-j|}$.
2. $\Omega_2$: AR(4), $\omega_{j'j} = \mathbf{I}(|j'-j|=0) + 0.4 \cdot \mathbf{I}(|j'-j|=1)$
   $+ 0.2 \cdot \mathbf{I}(|j'-j|=2) + 0.2 \cdot \mathbf{I}(|j'-j|=3) + 0.1 \cdot \mathbf{I}(|j'-j|=4)$.

3. $\Omega_3 = B + \delta I$, where each off-diagonal entry in $B$ is generated independently and equals 0.5 with probability $\alpha = 0.1$ or 0 with probability $1 - \alpha = 0.9$. $B$ has zeros on the diagonal, and $\delta$ is chosen so that the condition number of $\Omega_3$ is $p$ (keeping the diagonal constant across $p$ would result in either loss of positive definiteness or convergence to identity for larger $p$).
4. $\Omega_4$: Same as $\Omega_3$ except $\alpha = 0.5$.

All models are sparse (see Figure 2), and are numbered in order of decreasing sparsity (or increasing $s$). Note that the number of non-zero entries in $\Omega_1$ and $\Omega_2$ is proportional to $p$, whereas $\Omega_3$ and $\Omega_4$ have the expected number of non-zero entries proportional to $p^2$.

For all models, we generated $n = 100$ multivariate normal training observations and a separate set of 100 validation observations. We considered five different values of $p$, $30, 100, 200, 500$ and $1000$. The estimators were computed on the training data, with the tuning parameter for SPICE selected by minimizing the normal likelihood on the validation data. Using these values of the tuning parameters, we computed the estimated concentration matrix on the training data and compared it to the population concentration matrix.

We evaluate the concentration matrix estimation performance using the Kullback-Leibler loss,

$$\Delta_{KL}(\hat{\Omega}, \Omega) = \text{tr}\left(\Sigma\hat{\Omega}\right) - \log\left|\Sigma\hat{\Omega}\right| - p. \tag{28}$$

Note that this loss is based on $\hat{\Omega}$ and does not require inversion to compute $\hat{\Sigma}$, which is appropriate for a method estimating $\Omega$. The Kullback-Leibler loss was used by Yuan and Lin (2007) and Levina et al. (2008) to assess performance of methods estimating $\Omega$, and is obtained from the standard entropy loss of the covariance matrix (Lin and Perlman, 1985; Wu and Pourahmadi, 2003; Huang et al., 2006) by reversing the roles of $\Sigma$ and $\Omega$.

Results for the four covariance models are summarized in Table 1, which reports the average loss and the standard error over 50 replications. For $\Omega_1$, $\Omega_2$, and $\Omega_3$, SPICE outperforms the Ledoit-Wolf estimator for all values of $p$. The sample covariance performs much worse than either estimator in all cases (for $p = 30$). For $\Omega_4$, the least sparse of the four models, the Ledoit-Wolf estimator is about the same as SPICE (sometimes a little better, sometimes a little worse). This suggests, as we would expect from our bound on the rate of convergence, that SPICE provides the biggest gains in sparse models.

Table 1
*Simulations: Average (SE) Kullback-Leibler loss over 50 replications*

| $p$ | Sample | Ledoit-Wolf | SPICE | Sample | Ledoit-Wolf | SPICE |
|---|---|---|---|---|---|---|
| | | $\Omega_1$ | | | $\Omega_2$ | |
| 30 | 8.52(0.14) | 3.49(0.04) | 1.61(0.03) | 8.52(0.14) | 2.77(0.02) | 2.55(0.03) |
| 100 | NA | 26.65(0.08) | 8.83(0.05) | NA | 12.96(0.02) | 11.93(0.07) |
| 200 | NA | 76.83(0.13) | 21.23(0.09) | NA | 28.16(0.01) | 24.82(0.07) |
| 500 | NA | 262.8(0.19) | 78.26(0.26) | NA | 74.37(0.02) | 63.94(0.12) |
| 1000 | NA | 594.0(0.13) | 174.8(0.20) | NA | 151.9(0.04) | 133.7(0.20) |
| | | $\Omega_3$ | | | $\Omega_4$ | |
| 30 | 8.45(0.12) | 3.50(0.05) | 2.12(0.04) | 8.45(0.12) | 3.04(0.04) | 3.77(0.04) |
| 100 | NA | 29.25(0.44) | 17.09(0.10) | NA | 19.35(0.15) | 21.33(0.06) |
| 200 | NA | 86.93(1.64) | 45.58(0.13) | NA | 53.18(0.37) | 51.93(0.13) |
| 500 | NA | 240.3(3.24) | 168.7(0.37) | NA | 150.4(0.45) | 176.6(0.33) |
| 1000 | NA | 321.5(27.7) | 277.3(23.5) | NA | 269.8(18.1) | 307.3(20.6) |

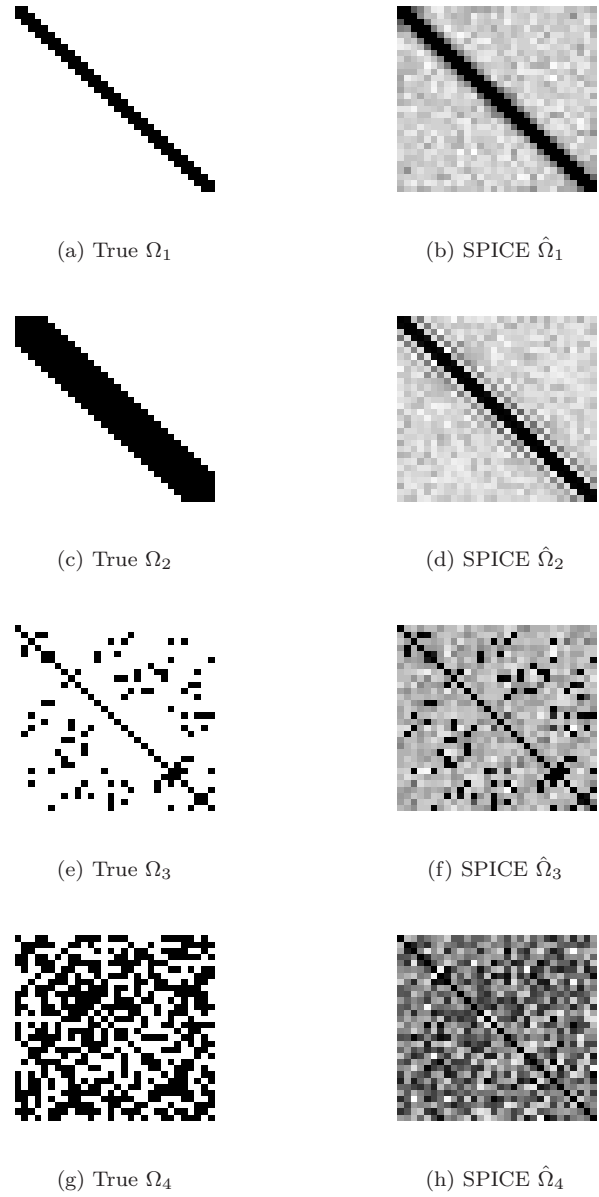

FIG 2. *Heatmaps of zeros identified in the concentration matrix out of 50 replications. White color is 50/50 zeros identified, black is 0/50.*

TABLE 2
*Percentage of correctly estimated non-zeros (TP %) and correctly estimated zeros (TN %) in the concentration matrix (average and SE over 50 replications) for SPICE*

| $p$ | TP % | TN % | TP % | TN % |
|---|---|---|---|---|
| | $\Omega_1$ | | $\Omega_2$ | |
| 30   | 100(0.00)   | 68.74(0.31) | 50.18(1.44) | 75.64(1.28) |
| 100  | 100(0.00)   | 74.70(0.08) | 49.96(1.10) | 72.68(1.21) |
| 200  | 100(0.00)   | 73.57(0.04) | 27.62(0.12) | 96.47(0.02) |
| 500  | 100(0.00)   | 91.97(0.01) | 22.48(0.09) | 98.81(0.00) |
| 1000 | 100(0.00)   | 98.95(0.00) | 22.29(0.05) | 98.82(0.00) |
| | $\Omega_3$ | | $\Omega_4$ | |
| 30   | 98.38(0.30) | 63.85(1.28) | 74.15(0.61) | 44.50(0.84) |
| 100  | 93.90(0.27) | 54.01(0.61) | 41.27(0.37) | 63.07(0.36) |
| 200  | 70.81(0.13) | 69.82(0.05) | 35.77(0.06) | 66.08(0.06) |
| 500  | 28.93(0.06) | 89.28(0.02) | 5.92(0.62)  | 94.27(0.61) |
| 1000 | 4.73(0.40)  | 72.36(6.13) | 2.07(0.14)  | 79.97(5.35) |

To assess the performance of SPICE on recovering the sparsity structure in the inverse, we report percentages of true non-zeros estimated as non-zero (TP %) and percentages of true zeros estimated as zero (TN %) in Table 2. We also plot heatmaps of the percentage of time each element was estimated as zero out of the 50 replications in Figure 2, for $p = 30$ for all four models. In general, recovering the sparsity structure is easier for smaller $p$ and for sparser models.

Finally, some example computing times: the SPICE algorithm for $\Omega_2$ takes about 2 seconds for $p = 200$, 1 minute for $p = 500$, and 15 minutes for $p = 1000$ on a regular PC. Glasso and SPICE both have complexity $O(p^3)$, but because of the quadratic approximation, SPICE tends to require more iterations to converge, and on average, we have observed a difference in computing times on the order of about 10 between glasso and SPICE. However, this factor does not grow with $p$, and SPICE computing times are still very reasonable even for large $p$.

## 4.2. Colon tumor classification example

In this section, we compare performance of covariance estimators for LDA classification of tumors using gene expression data from Alon et al. (1999). In this experiment, colon adenocarcinoma tissue samples were collected, 40 of which were tumor tissues and 22 non-tumor tissues. Tissue samples were analyzed using an Affymetrix oligonucleotide array. The data were processed, filtered, and reduced to a subset of 2,000 gene expression values with the largest minimal intensity over the 62 tissue samples. Additional information about the dataset and pre-processing can be found in Alon et al. (1999).

To assess the performance at different dimensions, we reduce the full dataset of 2,000 gene expression values by selecting $p$ most significant genes as mea-sured by the two-sample $t$-statistic, for $p = 50, 100, 200$. Then we use linear discriminant analysis (LDA) to classify these tissues as either tumorous or non-tumorous. We classify each test observation $\boldsymbol{x}$ to either class $k = 0$ or $k = 1$ using the LDA rule

$$\delta_k(\boldsymbol{x}) = \arg\max_k \left\{ \boldsymbol{x}^T \hat{\Omega} \hat{\boldsymbol{\mu}}_k - \frac{1}{2} \hat{\boldsymbol{\mu}}_k^T \hat{\Omega} \hat{\boldsymbol{\mu}}_k + \log \hat{\pi}_k \right\}, \tag{29}$$

where $\hat{\pi}_k$ is the proportion of class $k$ observations in the training data, $\hat{\boldsymbol{\mu}}_k$ is the sample mean for class $k$ on the training data, and $\hat{\Omega}$ is an estimator of the inverse of the common covariance matrix on the training data computed by one of the methods under consideration. Detailed information on LDA can be found in Mardia et al. (1979).

TABLE 3
*Averages and SEs of classification errors in % over 100 splits. Tuning parameter for SPICE chosen by (A): 5-fold CV on the training data maximizing the likelihood; (B): 5-fold CV on the training data minimizing the classification error; (C): minimizing the classification error on the test data*

|  |  | $p = 50$ | $p = 100$ | $p = 200$ |
|---|---|---|---|---|
| N. Bayes |  | 15.8(0.77) | 20.0(0.84) | 23.1(0.96) |
| L-W |  | 15.2(0.55) | 16.3(0.71) | 17.7(0.61) |
| SPICE | A | 12.1(0.65) | 18.7(0.84) | 18.3(0.66) |
| SPICE | B | 14.7(0.73) | 16.9(0.85) | 18.0(0.70) |
| SPICE | C | 9.0(0.57) | 9.1(0.51) | 10.2(0.52) |

To create training and test sets, we randomly split the data into a training set of size 42 and a testing set of size 20; following the approach used by Wang et al. (2007), we require the training set to have 27 tumor samples and 15 non-tumor samples. We repeat the split at random 100 times and measure the average classification error. The average errors with standard errors over the 100 splits are presented in Table 3. We omit the sample covariance because it is not invertible with such a small sample size, and include the naive Bayes classifier instead (where $\hat{\Sigma}$ is estimated by a diagonal matrix with sample variances on the diagonal). Naive Bayes has been shown to perform better than the sample covariance in high-dimensional settings (Bickel and Levina, 2004).

For an application such as classification, there are several possibilities for selecting the tuning parameter. Since we have no separate validation data available, we perform 5-fold cross-validation on the training data. One possibility (columns A in Table 3) is to continue using normal likelihood as a criterion for cross-validation, like we did in simulations. Another possibility (columns B in Table 3) is to use classification error as the cross-validation criterion, since that is the ultimate performance measure in this case. Table 3 shows that for SPICE both methods of tuning perform similarly. For reference, we also include the best error rate achievable on the test data, which is obtained by selecting the tuning parameter to minimize the classification error on the test data (columns C in Table 3). SPICE provides the best improvement over naive Bayes and Ledoit-Wolf for $p = 50$; for larger $p$, as less informative genes are added into the pool, the performance of all methods worsens.

## 5. Discussion

We have analyzed a penalized likelihood approach to estimating a sparse concentration matrix via a lasso-type penalty, and showed that its rate of convergence depends explicitly on how sparse the true matrix is. This is analogous to results for banding (Bickel and Levina, 2008), where the rate of convergence depends on how quickly the off-diagonal elements of the true covariance decay, and for thresholding (Bickel and Levina, 2007; El Karoui, 2007), where the rate also depends on how sparse the true covariance is by various definitions of sparsity. We conjecture that other structures can be similarly dealt with, and other types of penalties may show similar behavior when applied to the "right" type of structure – for example, a ridge, bridge, or other more complex penalty may work well for a model that is not truly sparse but has many small entries. A generalization of this work to other penalties has been recently completed by Lam and Fan (2007), who have also proved "sparsistency" of SPICE-type estimators.

While we assumed normality, it can be replaced by a tail condition, analogously to Bickel and Levina (2008). The use of normal likelihood is, of course, less justifiable if we do not assume normality, but it was found empirically that it still works reasonably well as a loss function even if the true distribution is not normal (Levina et al., 2008).

The Cholesky decomposition of covariance was only considered appropriate when variables are ordered, and we have shown it to be a useful tool for enforcing positive definiteness of the estimator even when variables have no natural ordering. Our optimization algorithm has complexity of $O(p^3)$ and is equally applicable to general $l_q$ penalties.

## Acknowledgments

We thank an Associate Editor and two referees for their feedback and helpful comments, Sourav Chatterjee (Berkeley) for helpful discussions, Noureddine El Karoui (Berkeley) for helpful discussions and corrections, and Shuheng Zhou (CMU) for a correction. P. J. Bickel's research is partially supported by a grant from the NSF (DMS-0605236). E. Levina's research is partially supported by grants from the NSF (DMS-0505424 and DMS-0805798) and the NSA (MSPF-04Y-120). J. Zhu's research is partially supported by grants from the NSF (DMS-0505432 and DMS-0705532).

## Appendix A: Derivation of the Algorithm

In this section we give a full derivation of the parameter update equations involved in the optimization algorithm. Recall that we have re-parametrized the objective function (20) using (22)–(24). We cycle through the parameters in $T$ and for each $t_{lc}$, compute partial derivatives with respect to $t_{lc}$ while holding all other parameters fixed, and solve the univariate linear equation corresponding to setting this partial derivative to 0.

For simplicity, we separate the likelihood and the penalty by writing $\tilde{f}(T) = \ell(T) + P(T)$. We also suppress the $\epsilon$-perturbation in the denominator for simplicity of notation. For the likelihood part, taking the partial derivative with respect to $t_{lc}$, $1 \leq c \leq p$, $c \leq l \leq p$ gives

$$\frac{\partial}{\partial t_{\ell c}}\ell(T) = -2\underbrace{\frac{\partial}{\partial t_{\ell c}}\sum_{j=1}^{p}\log t_{jj}}_{=0 \text{ if } j \neq c} + \frac{1}{n}\sum_{i=1}^{n}\underbrace{\frac{\partial}{\partial t_{\ell c}}\sum_{j=1}^{p}\left(\sum_{k=1}^{j}t_{jk}X_{ik}\right)^2}_{=0 \text{ if } j \neq l}$$

$$= \frac{-2}{t_{cc}}\mathbf{I}\{l = c\} + t_{lc}[2\hat{\sigma}_{cc}] + 2\sum_{k=1,\ k\neq c}^{l}t_{lk}\hat{\sigma}_{kc}, \qquad (30)$$

For the penalty part, using the quadratic approximation (25) gives

$$\frac{\partial}{\partial t_{\ell c}}P(T) \approx \frac{\partial}{\partial t_{\ell c}}\sum_{j'>j}\frac{\lambda q}{|\omega_{j'j}^0|^{2-q}}\omega_{j'j}^2 = \sum_{k=1,k\neq c}^{l}\frac{\lambda q}{|\omega_{ck}^0|^{2-q}}\frac{\partial}{\partial t_{\ell c}}\omega_{ck}^2, \qquad (31)$$

since the only nonzero terms in (31) are those for which $j' \leq l$ and either $j' = c$ or $j = c$. For $1 \leq k \leq l$ such that $k \neq c$, we have $\frac{\partial}{\partial t_{\ell c}}\omega_{ck}^2 = 2\omega_{ck}t_{lk}$, and collecting terms together we get

$$\frac{\partial}{\partial t_{\ell c}}P(T) = t_{lc}\left[2\lambda q\sum_{k=1,k\neq c}^{l}\frac{t_{lk}^2}{|\omega_{ck}^0|^{2-q}}\right] + 2\lambda q\sum_{k=1,k\neq c}^{l}\frac{(\omega_{ck}-t_{lc}t_{lk})t_{lk}}{|\omega_{ck}^0|^{2-q}}. \qquad (32)$$

Combining together (30) and (32), we have the parameter update equation for $t_{lc}$ when $l \neq c$, is given by

$$\hat{t}_{lc} = \frac{-\sum_{k=1,\ k\neq c}^{l}t_{lk}\hat{\sigma}_{kc} - \lambda q\sum_{k=1,k\neq c}^{l}(\omega_{ck}-t_{lc}t_{lk})t_{lk}|\omega_{ck}^0|^{q-2}}{\hat{\sigma}_{cc} + \lambda q\sum_{k=1,k\neq c}^{l}t_{lk}^2|\omega_{ck}^0|^{q-2}}.$$

If $l = c$, we solve $au^2 + bu - 1 = 0$ for $u$ using the quadratic formula, where

$$a = \hat{\sigma}_{cc} + \lambda q\sum_{k=1,k\neq c}^{l}t_{lk}^2|\omega_{ck}^0|^{q-2},$$

$$b = \sum_{k=1,\ k\neq c}^{l}t_{lk}\hat{\sigma}_{kc} + \lambda q\sum_{k=1,k\neq c}^{l}(\omega_{ck}-t_{lc}t_{lk})t_{lk}|\omega_{ck}^0|^{q-2},$$

then take the positive solution $\hat{t}_{cc} = u^+$.

We also can quickly update the $\omega_{ck}$ involving $t_{lc}$ via

$$\omega_{ck} = \omega_{ck}^0 + t_{lk}(\hat{t}_{lc} - t_{lc}).$$